\documentclass{amsart}
\usepackage{amssymb, amsmath,color}

\newtheorem{theorem}{Theorem}
\newtheorem{lemma}[theorem]{Lemma}

\newtheorem{remark}[theorem]{Remark}

\begin{document}

\title{Locally Conformally Flat Lorentzian Gradient Ricci Solitons}
\author{M. Brozos-V\'{a}zquez $\,$ E. Garc\'{i}a-R\'{i}o $\,$ S. Gavino-Fern\'{a}ndez}
\address{MBV: Department of Mathematics, University of A Coru\~na, Spain}
\email{miguel.brozos.vazquez@udc.es}
\address{EGR, SGF: Faculty of Mathematics,
University of Santiago de Compostela,
15782 Santiago de Compostela, Spain}
\email{eduardo.garcia.rio@usc.es \,\, sandra.gavino@usc.es}
\thanks{2010 {\it Mathematics Subject Classification}: 53C21, 53C50, 53C25.\\
M. B.-V. and E. G.-R. are supported by projects MTM2009-07756 and INCITE09 207 151 PR (Spain). S. G.-F. is supported by project MTM2009-14464-C02-01 (Spain)}
\keywords{Ricci solitons, gradient Ricci solitons, Lorentzian locally conformally flat manifolds}

\begin{abstract}
It is shown that locally conformally flat Lorentzian gradient Ricci solitons are locally isometric to a Robertson-Walker warped product, if the gradient of the potential function is non null, and to a plane wave, if the gradient of the potential function is null. The latter gradient Ricci solitons are necessarily steady.
\end{abstract}

\maketitle

\section{Introduction}\label{se:1}

Let $M$ be a differentiable manifold of dimension $n+2$, let $g$ be a pseudo-Riemannian metric and let $f$ be a { smooth} function on $M$. We say that the triple $(M,g,f)$ is a \emph{gradient Ricci soliton} if the following equation is satisfied:
\begin{equation}\label{gradsoliton}
\text{Hes}_{f}+\rho=\lambda g,
\end{equation}
where $\operatorname{Hes}$ denotes de Hessian, $\rho$ denotes the Ricci tensor and $\lambda$ is a real number. By contracting in the previous equation one sees that $\lambda=\frac{1}{n+2}(\Delta f+\tau)$, where $\tau$ denotes the scalar curvature and $\Delta$ denotes the Laplacian. A gradient Ricci soliton is said to be \emph{shrinking, steady} or \emph{expanding} if $\lambda >0,$ $\lambda =0$ or $\lambda <0$, respectively.

Gradient Ricci solitons are a particularly interesting family of Ricci solitons. These arise as self-similar solutions of the Ricci flow
$\frac{\partial}{\partial t}g(t)=-2\rho_{g(t)}$ under certain conditions.
Lorentzian Ricci solitons have been investigated recently
showing many essential differences with respect to the positive definite case \cite{BBGG, BCGG, onda}.

A gradient Ricci soliton $(M,g,f)$ is said to be \emph{rigid} if $(M,g)$ is isometric to a quotient of $N\times\mathbb{R}^k$, where $N$ is an Einstein manifold and the potential function $f$ is defined on the Euclidean factor as $f=\frac{\lambda}{2}\| x\|^2$ (thus generalizing the Gaussian soliton) \cite{PW3}. Although rigidity is a rather restrictive condition, rigid Ricci solitons are the only solitons in many important situations as shown in \cite{PW2}, where it is proven that any homogeneous gradient Ricci soliton is rigid if the metric is positive definite.
This result fails when passing from the Riemannian to the Lorentzian setting \cite{BCGG}. Indeed, indecomposable Lorentzian symmetric spaces provide examples of nontrivial steady gradient Ricci solitons in which, moreover, the gradient of the potential function $\nabla f$ is a null geodesic vector field \cite{BBGG}.

Riemannian locally conformally flat complete shrinking and steady gradient Ricci solitons were recently classified: they are quotients of $\mathbb{R}^{n+2}$, $\mathbb{S}^{n+2}$ or $\mathbb{R}\times\mathbb{S}^{n+1}$ if shrinking and the Bryant soliton if steady \cite{Cao-Chen} (see also \cite{FL-GR}). The existence of locally conformally flat Lorentzian steady gradient Ricci solitons of non Bryant type was proven in \cite{BBGG}.

The purpose of this work is to investigate locally conformally flat gradient Ricci solitons in the Lorentzian setting by focusing on their local structure. The following is the main result.

\begin{theorem}\label{mainth}
Let $(M,g,f)$ be a locally conformally flat Lorentzian gradient Ricci soliton.
\begin{enumerate}
\item[(i)] In a neighborhood of any point where $\|\nabla f\|\neq 0$, $M$ is locally isometric to a
Robertson-Walker warped product $I\times_\psi N$ with metric $\varepsilon dt^2+\psi^2 g_N$, where $I$ is a real interval and $(N,g_N)$ is a space of constant curvature $c$.

\item[(ii)] If $\|\nabla f\|=0$ on a non-empty open set, then $(M,g)$ is locally isometric to a plane wave, i.e., $M$ is locally diffeomorphic to $\mathbb{R}^2\times\mathbb{R}^n$ with metric
    \[
    g=2dudv +H(u,x_1,\dots,x_n) du^2 +\sum _{i=1}^n dx _i ^2,
    \]
    where $H(u,x_1,\dots,x_n)= a(u) \sum _{i=1}^n x_i ^2 + \sum _{i=1}^n b_i (u) x_i +c(u)$
    for some functions $a(u)$, $b_i(u)$, $c(u)$
and the potential function is given by $f(u,x_1,\dots,x_n)$ $=$ $f_0(u)$, with $f_0''(u)=-\rho_{uu} =n\,a(u)$.
\end{enumerate}
\end{theorem}

We say that a gradient Ricci soliton is \emph{non isotropic} if $\|\nabla f\|\neq 0$ or
\emph{isotropic} if $\|\nabla f\|=0$. In the following sections we will study both cases separately.

\begin{remark}\rm
Riemannian locally conformally flat gradient Ricci solitons are analogous to the manifolds describe in Theorem~\ref{mainth}-(i). Due to holonomy action there exist other possibilities in Lorentzian signature, as Theorem~\ref{mainth}-(ii) shows.
\end{remark}

\begin{remark}\rm
The character of $\nabla f$ may vary from one point to another. Thus, for example, consider the { Lorentzian analog of the Gaussian soliton. Let $(\mathbb{L}^{n+2},g)$ be the flat Minkowski space and let $f(x_1,\dots,x_{n+2})=\frac{\lambda}2(-x_1^2+x_2^2+\dots+x_{n+2}^2)$ be defined on $\mathbb{L}^{n+2}$}. The gradient of $f$ is given by $\nabla f=\lambda(x_1+x_2+\dots+x_{n+2})$ and the Hessian is $Hes_f=\lambda g$. Hence the soliton equation \eqref{gradsoliton} is satisfied for any given $\lambda$. Note that $\|\nabla f\|^2=\lambda^2 (-x_1^2+x_2^2+\dots+x_{n+2}^2)$ is positive, zero or negative depending on $(x_1,\dots,x_{n+2})$, so the character of $\nabla f$ varies with the point.
\end{remark}

The paper is organized as follows. In Section~\ref{se:2} we recall some basic formulas and give some results showing that under certain assumptions $\nabla f$ is an eigenvector of the Ricci operator; this will be crucial in the proof of Theorem~\ref{mainth}. Also we study two-dimensional gradient Ricci solitons and Einstein gradient Ricci solitons. We devote Section~\ref{se:3} to analyze locally conformally flat non isotropic gradient Ricci solitons and Section~\ref{se:4} to study the isotropic { case, showing that the underlying structure of such a soliton is a $pp$-wave. Finally the existence of gradient Ricci solitons in $pp$-waves is discussed in general, without any further assumption, in Section~\ref{se:5}. Thus the restriction of this discussion} to locally conformally flat $pp$-waves completes the proof of Theorem~\ref{mainth}.

\section{General formulae and remarks}\label{se:2}

The orthogonal group decomposes the space of curvature tensors into three irreducible modules. Thus, a curvature tensor $R$ can be written as $R=\frac{\tau}{2(n+2)(n+1)}g\odot g + \frac{1}{n}\rho_0\odot g +W$, where $\rho_0$ is the traceless Ricci tensor, $W$ is the Weyl conformal tensor, and $\odot$ is the Kulkarni-Nomizu product (for $A$, $B$ symmetric $2$-tensors,  $(A\odot B)(x,y,z,w)=A(x,z)B(y,w)+A(y,w)B(x,z)-A(x,w)B(y,z)-A(y,z)B(x,w)$). Note that the curvature tensor can also be written as $R=C\odot g+W$ where $C=\frac{1}{n}\left(\rho-\frac{\tau}{2(n+1)}g\right)$ is the Schouten tensor. The summands in the previous decomposition have a geometrical meaning; thus, for example, Einstein manifolds have vanishing traceless Ricci tensor, while locally conformally flat manifolds have $W=0$ if $n\geq 2$ and the Schouten tensor is Codazzi (i.e., its covariant derivative is totally symmetric) if $n\geq 1$. In this Section we begin the study of gradient Ricci solitons in these two
  particular cases.

Let $(M,g,f)$ be a Lorentzian gradient Ricci soliton. Although the following is well-known (see, for example, \cite{PW3}), we recall the formulae and sketch the proof in order to make the paper as self-contained as possible. Let $\operatorname{Ric}$ denote the Ricci operator defined by $g(\operatorname{Ric}(X),Y)=\rho(X,Y)$ for any vector fields $X$ and $Y$.

\begin{lemma}\label{lemma:formulas}
A Lorentzian gradient Ricci soliton with potential function $f$ satisfies
\begin{align}
\nabla \tau =2 \operatorname{Ric}(\nabla f),\label{condgradsol1}\\
\tau+\|\nabla\, f\|^2-2\lambda f=\operatorname{const}.\label{condgradsol2}
\end{align}
\end{lemma}

\begin{proof}
Tracing the soliton equation (\ref{gradsoliton}) gives $\Delta f+\tau=n\lambda$,  hence
$\nabla \tau=-\nabla \Delta f$.
The contracted second Bianchi identity ($\nabla_Z \tau=2\text{div} \rho(Z)$) together with the Bochner formula $\text{div}(\nabla \nabla f)=\rho(\nabla f)+\nabla \Delta f$ now gives \eqref{condgradsol1}.

Writing the soliton equation as $\operatorname{Ric}+\nabla \nabla f=\lambda Id$ and using (\ref{condgradsol1}) one has
\[
\nabla \tau=2 \operatorname{Ric}(\nabla f)=2\lambda \nabla f-2 \nabla_{\nabla f}\nabla f=2\lambda \nabla f-\nabla \|\nabla f\|^2.
\]
Hence $\nabla(\tau+\|\nabla f\|^2-2\lambda f)=0$, which proves \eqref{condgradsol2}.
\end{proof}

\begin{remark}\rm\label{remark:eigenvectors}
As a consequence of Lemma~\ref{lemma:formulas}, there are several particular situations in which $\nabla f$ is an eigenvector of the Ricci operator. Thus, if $\tau$ is constant, from \eqref{condgradsol1} it follows that $\nabla f$ is an eigenvector for the Ricci operator associated to the eigenvalue zero. Also, if $\nabla f$ is null, then from \eqref{condgradsol2} one has $\tau=\operatorname{const} + 2\lambda f$; now substitute in \eqref{condgradsol1} to see that $\operatorname{Ric}(\nabla f)=\lambda \nabla f$.
\end{remark}

\subsection{Two-dimensional steady gradient Ricci solitons}

Let $(M,g,f)$ be a two-dimen\-sio\-nal gradient Ricci soliton. Consider the canonical para-K\"ahler structure $\mathfrak{J}$ (i.e., $\mathfrak{J}^2=Id$, $g(\mathfrak{J}\,\cdot\,,\mathfrak{J}\,\cdot\,)=-g(\,\cdot\,,\,\cdot\,)$, $\nabla\mathfrak{J}=0$) on $(M,g)$. Then, proceeding as in \cite[$\S 1.3$]{CCGGIIKLLN}, it follows that $\mathfrak{J}\nabla f$ is a Killing vector field.
Now, if $\nabla f$ is a nonnull vector field, then $(M,g)$ is locally a warped product. On the other hand, if $\nabla f$ is a null vector field, then consider coordinates $(x_1,x_2)$ so that the null Killing vector field is
$\mathfrak{J}\nabla f=\frac{\partial}{\partial x_2}$. The metric tensor takes the form
$g=a(x_1,x_2)dx_1^2+b(x_1,x_2)dx_1 dx_2$ for some functions $a$, $b$. The fact that $\frac{\partial}{\partial x_2}$ is Killing implies that  $\frac{\partial}{\partial x_2}a=0$ and
$\frac{\partial}{\partial x_2}b=0$, therefore $g$ is flat and $\mathfrak{J}\nabla f$ is indeed parallel.

Next we are going to extend the Hamilton cigar soliton to the Lorentzian setting. Let $(M,g,f)$ be a two-dimensional steady gradient Ricci soliton with $\nabla f$ a timelike vector field (the spacelike case is similar). Set $M=I\times N$ with metric $g=-dt^2+\omega(t)^2ds^2$ and assume $f$ only depends on $t$.
Then a straightforward calculation from \eqref{gradsoliton} gives that $(M,g,f)$ is a steady gradient Ricci soliton if and only if
\[
f''(t)-\frac{\omega''(t)}{\omega(t)}=0
\quad\text{and}\quad
-f'(t) \omega'(t)+\omega''(t)=0.
\]
Hence $f''w-f'w'=0$, so we integrate to see that $f'(t)=\kappa w(t)$ { for a constant $\kappa$}. Equations above reduce to $\kappa \omega\omega'-\omega''=0$.
Hence the possible solutions depending on the sign of $\kappa$ are given by
\begin{itemize}
\item[(i)] If $\kappa=0$; then $\omega(t)=a t+b$ for constants $a$ and $b$. In this case $M$ is flat { and $f$ is constant}.
\item[(ii)] If $\kappa=r^2$; then $\omega(t)=\frac{a\sqrt{2}}{r}\tan{r\sqrt{2} (a t+b)}$ where $a$ and $b$ are constants. The potential function is
$f(t)=d-2\log{\cos(\frac{r (a t+b)}{\sqrt{2}})}$ for a constant $d$, and the scalar curvature is $\tau=2a\,r^2 \,\text{sec}^2(r \frac{(at+b)}{\sqrt{2}})$.
\item[(iii)] If $\kappa=-r^2$; then $\omega(t)=\frac{a\sqrt{2}}{r}\tanh{r\sqrt{2} (a t+b)}$ for $a$ and $b$ constants. The potential function is
$f(t)=d+2\log{\cosh(\frac{r (a t+b)}{\sqrt{2}})}$ for a constant $d$, and the scalar curvature is $\tau=-2a\,r^2 \,\text{sech}^2(r \frac{(at+b)}{\sqrt{2}})$.
\end{itemize}

Analyzing geodesic completeness in the Lorentzian case is a subtle task. Indeed Lorentzian warped products of geodesically complete manifolds need not be complete, { as occurs in positive definite signature}. Necessary and sufficient conditions for geodesic completeness of Lorentzian warped products were investigated in \cite{CS}. Let $M=I\times_\omega N$ be a warped product where $I=(\alpha,\beta)$ is a real interval and $(N,g_N)$ is a geodesically complete manifold. Then $M$ is timelike, spacelike and null geodesically complete if and only if for some $\gamma\in (\alpha,\beta)$ it follows that
\[
\int_\alpha^\gamma \frac{\omega}{\sqrt{1+\omega^2}}dt=
\int_\gamma^\beta \frac{\omega}{\sqrt{1+\omega^2}}dt=+\infty.
\]
As a consequence, Lorentzian warped products given by (ii) above are not geodesically complete, while those given by (iii) are. Thus, (iii) generalizes Hamilton's cigar (see \cite{hamilton}) to the Lorentzian setting.

\subsection{Lorentzian Einstein gradient Ricci solitons}

Ricci solitons are generalizations of Einstein metrics. If $(M,g)$ is a complete Riemannian Einstein manifold, then $(M,g,f)$ is a  gradient Ricci soliton if and only if it has $\operatorname{Hes}_f=0$ or is a Gaussian \cite{PW3}. The next result describes the local structure of Einstein gradient Ricci solitons in the Lorentzian setting.

\begin{theorem}\label{th:Einstein}
Let $(M,g)$ be a Lorentzian Einstein manifold. If $(M,g,f)$ is a gradient Ricci soliton with nonconstant $f$, then $(M,g)$ is Ricci flat. Moreover:
\begin{itemize}
\item[(i)] If  $\|\nabla f\|\neq 0$, then $(M,g)$ is locally a warped product of the form $I\times_{f'}N$  and
    the potential function $f(t)=\frac{\lambda}{2}t^2+at+b$.
\item[(ii)] If $\|\nabla f\|=0$, then
there exist coordinates $(u,v,x_1,\dots,x_n)$ in which the metric has the form
$g=2dudv +\tilde g$, where the $n$-dimensional metric $\tilde g$ does not depend on $v$. Moreover, the potential function $f$ is given by any function $f(u)$ with $f''(u)=0$ and the soliton is steady.
\end{itemize}
\end{theorem}

\begin{proof}
Let $(M,g,f)$ be an Einstein gradient Ricci soliton. Then it follows from \eqref{condgradsol1} that either the potential function $f$ is constant or otherwise $(M,g)$ is Ricci flat. Assume $(M,g)$ is Ricci flat. The soliton equation
\eqref{gradsoliton} reduces to $\text{Hes}_{f}=\lambda g=\frac{\Delta f}{n+2}g$. This equation was previously investigated by Brinkmann \cite{brinkmann} (see \cite{KR} for a modern exposition) showing that in a neighborhood of any point where $\|\nabla f\|\neq 0$ the manifold $(M,g)$ decomposes locally as a warped product of a real interval $I\subset\mathbb{R}$ and an Einstein manifold $(N,g_N)$ so that
$g=\varepsilon dt^2+(f')^2g_N$, where $f$ is a real function defined on $I$ with $f'\neq 0$. Now, since $(M,g)$ is Ricci flat, a direct computation of the Ricci tensor for the metric $\varepsilon dt^2+(f')^2g_N$ shows that $g_N$ is Einstein and $f$ must satisfy $f'''=0$ and $f'f'''+n\varepsilon(f'')^2=\frac{\tau_N}{n+1}$. Hence $f(t)=\frac{\lambda}{2}t^2+at+b$ and $\tau_N=n(n+1)\varepsilon \lambda^2$.

Now assume $\|\nabla f\|=0$ identically. Then \eqref{condgradsol2} shows that either $f$ is constant or the gradient Ricci soliton is steady. If $\lambda=0$ the Ricci soliton equation \eqref{gradsoliton} reduces to $\operatorname{Hes}_{f}=0$. Then $\nabla f$ is a parallel isotropic vector field and the metric tensor can be written in suitable Rosen coordinates $(u,v,x_1,\dots,x_n)$ as $g=dudv+\tilde g$, where the $n$-dimensional metric $\tilde g(u)$ is Ricci flat for any fixed $u$ and does not depend on $v$ \cite{brinkmann, KR}. Moreover, in this coordinates $\nabla f=\frac{\partial}{\partial v}$ and the potential function depends only on the variable $u$. Now, the result follows by computing the Hessian of $f$.
\end{proof}

The (not complete) gradient Ricci solitons described in Theorem~\ref{th:Einstein}-(i) do have a Riemannian analog. However, those given in Theorem~\ref{th:Einstein}-(ii) are steady gradient Ricci solitons $(M,g,f)$, with $\operatorname{Hes}_f=0$, without Riemannian counterpart.

\subsection{General remarks on locally conformally flat gradient Ricci solitons}
Although locally conformally flat gradient Ricci solitons will be more deeply analyzed in Sections~\ref{se:3} and \ref{se:4}, we begin here by establishing a technical lemma. Proceeding in a similar way to { that developed} in \cite{FL-GR}, one has the following:

\begin{lemma}\label{le:5}
Let $(M,g,f)$ be a locally conformally flat gradient Ricci soliton. Then $\nabla f$ is an eigenvector of the Ricci operator.
\end{lemma}

\begin{proof}
Since $(M,g)$ is locally conformally flat the Schouten tensor is Codazzi, this is, $(\nabla_X C)(Y,Z)=(\nabla_Y C)(X,Z)$ for all vector fields $X,Y,Z$.
Hence
\begin{align}\label{codazzi1}
(\nabla_X \rho)(Y,Z)-\frac{X(\tau)}{2(n+1)}g(Y,Z)=(\nabla_Y \rho)(X,Z)-\frac{Y(\tau)}{2(n+1)}g(X,Z).
\end{align}

From (\ref{gradsoliton}) and using that $\text{Hes}_f(X,Y)=g(\nabla_X \nabla f,Y)$ one has
\begin{align*}
(\nabla_X \rho)(Y,Z)& =-(\nabla_X \text{Hes}_f)(Y,Z)\\
  & =-X g(\nabla_Y \nabla f,Z)+g(\nabla_{\nabla_X Y}\nabla f,Z)+g(\nabla_Y \nabla f,\nabla_X Z)\\
  & =-g(\nabla_X \nabla_Y \nabla f,Z)+g(\nabla_{\nabla_X Y}\nabla f,Z).
\end{align*}

Substituting this expression in  (\ref{codazzi1}) we get
\begin{align*}
g(\nabla_X \nabla_Y \nabla f,Z)&-g(\nabla_{\nabla_X Y}\nabla f,Z)+\frac{X(\tau)}{2(n+1)}g(Y,Z)\\
 & =g(\nabla_Y \nabla_X \nabla f,Z)-g(\nabla_{\nabla_Y X}\nabla f,Z)+\frac{Y(\tau)}{2(n+1)}g(X,Z).
\end{align*}
Thus
\[
g(\nabla_X\nabla_Y\nabla f-\nabla_Y\nabla_X\nabla f-\nabla_{[X,Y]}\nabla f,Z)=-\frac{X(\tau)}{2(n+1)}g(Y,Z)+\frac{Y(\tau)}{2(n+1)}g(X,Z),
\]
that is,
\[
R(X,Y,Z,\nabla f)=-\frac{X(\tau)}{2(n+1)}g(Y,Z)+\frac{Y(\tau)}{2(n+1)}g(X,Z),
\]
or equivalently, using \eqref{condgradsol1},
\begin{equation}\label{curv-simplify}
R(X,Y,Z,\nabla f)=-\frac{1}{n+1}\rho(X,\nabla f)g(Y,Z)+\frac{1}{n+1}\rho(Y,\nabla f)g(X,Z).
\end{equation}
Let $Z=\nabla f$ in (\ref{curv-simplify}) to obtain
\[
\rho(Y,\nabla f)g(X,\nabla f)=\rho(X,\nabla f)g(Y,\nabla f){ .}
\]
Now choose $X$ so that $g(X,\nabla f)=1$ to see that for all $Y\perp \nabla f$ one has
\[
0=\rho(Y,\nabla f)=-\operatorname{Hes}_f(Y,\nabla f)
\]
and conclude that $\nabla f$ is an eigenvector of the Ricci operator. Note that $\nabla f$ is also an eigenvector of the Hessian operator $\text{hes}_f(X)=\nabla_X\nabla f$.
\end{proof}

\section{Non isotropic gradient Ricci solitons}\label{se:3}

Next we show that in
a neighborhood of any point where  $\|\nabla f\|\neq 0$ the underlying manifold has the local structure of a warped product,
thus proving Theorem~{\ref{mainth}~-~$(i)$}.

\begin{lemma}\label{lemma:non-isotropic}
Let $(M,g,f)$ be a locally conformally flat Lorentzian gradient Ricci soliton with $\|\nabla f\|_P\neq 0$ for some point $P\in M$. Then, on a neighborhood of $P$, $(M,g)$ is a warped product of a real interval and a space of constant sectional curvature $c$.
\end{lemma}

\begin{proof}
If the Weyl tensor of $(M,g)$ vanishes, then the curvature tensor expresses as
\begin{equation}\label{curv-general}
\begin{array}{rcl}
R(X,Y,Z,T)& =& \displaystyle\frac{\tau}{n(n+1)}\left\{g(X,T)g(Y,Z)-g(X,Z)g(Y,T)\right\}\\
\noalign{\medskip}
            &&\phantom{\frac{\tau}{n(n+1)}\left\{\right.}+\frac{1}{n}\left\{\rho(X,Z)g(Y,T)+\rho(Y,T)g(X,Z)\right.\\
\noalign{\medskip}
            &&\phantom{\frac{\tau}{n(n+1)}+\frac{1}{n}\left\{\right.}\left.-\rho(X,T)g(Y,Z)-\rho(Y,Z)g(X,T)\right\}.
\end{array}
\end{equation}

Consider the unit vector $V=\frac{\nabla f}{\|\nabla f\|}$, which can be timelike or spacelike (we set $g(V,V)=\varepsilon$), on the tangent space $T_PM$. Complete it to a local orthonormal frame
$\{V,E_1,\dots,E_{n+1}\}$ with $g(E_i,E_i)=\varepsilon_i$. Then  from (\ref{curv-simplify}) one has
\[
R(V,E_i,E_i,V)=-\frac{1}{n+1}\rho(V,V)\varepsilon_i\,,
\]
{ while} from (\ref{curv-general}) one gets
\[
R(V,E_i,E_i,V)=\frac{\tau}{n(n+1)}\varepsilon\varepsilon_i -\frac{1}{n}\rho(V,V)\varepsilon_i-\frac{1}{n}\rho(E_i,E_i)\varepsilon\,.
\]

\noindent Hence for all $i=1,\dots,n+1$:
\[
-\frac{1}{n+1}\rho(V,V)\varepsilon_i=-\frac{1}{n}\rho(V,V)\varepsilon_i-\frac{1}{n}\rho(E_i,E_i)\varepsilon +\frac{\tau}{n(n+1)}\varepsilon\varepsilon_i,
\]
from where $\rho(E_i,E_i)\varepsilon=\frac{1}{n+1}(\tau \varepsilon-\rho(V,V))\varepsilon_i${ . Using} \eqref{gradsoliton} we have
\[
\operatorname{Hes}_f(E_i,E_i)=\lambda \varepsilon_i+\frac{1}{n+1}\left(\rho(V,V)\varepsilon-\tau \right)\varepsilon_i,
\]
which shows that the level sets of $f$ are totally umbilical hypersurfaces. Hence $(M,g)$ decomposes locally as a twisted product of the form $I\times_\omega N$
(see  \cite[Thm. 1]{PR}). Now, since $\nabla f$ is an eigenvector of the Ricci operator by Lemma~\ref{le:5}, it follows that $\rho(V,E_i)=0$ for all $i=1,\dots,n+1$,
and therefore the twisted product reduces to a warped product \cite[Thm. 1]{FGKU}. Hence $(M,g)$ is locally a warped product
$(I\times N,\varepsilon dt^2+\psi(t)^2 g_N)$ where $(N,g_N)$ is a Riemannian or a Lorentzian manifold of constant sectional curvature $c$ \cite{MER}.
\end{proof}

\begin{remark}\rm
The potential function $f$ in Lemma~\ref{lemma:non-isotropic} is a radial function $f(t)$, and hence a direct computation from the soliton equation \eqref{gradsoliton} shows that it is given as a solution to the equations:
\[
f''=\varepsilon\lambda + (n+1)\frac{\psi}{\psi''}, \qquad
 \varepsilon\psi\psi' f'=\lambda \psi^2 -nc+\varepsilon (\psi\psi''+n(\psi')^2).
\]
Note that this equations impose restrictions on the warping function $\psi$, thus the warped product is not arbitrary.
\end{remark}

\section{Locally conformally flat isotropic gradient Ricci solitons}\label{se:4}

{ In this section we consider the case of gradient Ricci solitons with $\|\nabla f\|=0$.}

Recall that for a Riemannian metric the holonomy group acts completely reducibly, i.e., the tangent space decomposes into subspaces on which it acts trivially or irreducibly, but for indefinite metrics the situation is more subtle. Indecomposable but not irreducible Lorentzian manifolds admit a parallel degenerate line field $\mathcal{D}$, and thus
the curvature satisfies (see, for example, \cite{DR})
\begin{equation}
R(\mathcal{D},\mathcal{D}^\perp,\cdot,\cdot)=0,\quad
R(\mathcal{D},\mathcal{D},\cdot,\cdot)=0, \quad \text{ and } \, R(\mathcal{D}^\perp,\mathcal{D}^\perp,\mathcal{D},\cdot)=0.
\end{equation}
If $\mathcal{D}$ is spanned by a parallel null vector field, then $(M,g)$ is said to be a \emph{$pp$-wave} if
\begin{equation}\label{leistner}
R(\mathcal{D}^\perp,\mathcal{D}^\perp,\cdot,\cdot)=0.
\end{equation}
$(M,g)$ is called a \emph{$pr$-wave} if \eqref{leistner} is satisfied but $\mathcal{D}$ is not necessarily spanned by a parallel vector field. Clearly any $pp$-wave is a $pr$-wave, and the converse is true if the Ricci tensor is isotropic (i.e., the image of the Ricci operator is totally isotropic) \cite{leistner}.
The general form of an $(n+2)$-dimensional $pp$-wave is the following: the ambient space is
$\mathbb{R}^{n+2}$ with coordinates $(u,v,x_1,..,x_{n})$, $n \geq 1$, and the Lorentzian metric is given by
\begin{equation}\label{pp}
g_{ppw}=2dudv +H(u,x_1,\dots,x_{n}) du^2 +\sum _{i=1}^n dx _i ^2,
\end{equation}
where $H(u,x_1,..,x_{n})$ is an arbitrary smooth function.

\begin{lemma}\label{lemma:isotropic}
Any isotropic locally conformally flat Lorentzian gradient Ricci soliton is steady and the underlying manifold is locally a $pp$-wave.
\end{lemma}

\begin{proof}
Let $(M,g,f)$ be a gradient Ricci soliton with $\|\nabla f\|=0$. In what follows we will show that $\nabla f$ spans a parallel null line field and { furthermore} that \eqref{leistner} holds.
Set $V=\nabla f$.
Since $V$ is a null vector, there exist orthogonal vectors $S$, $T$ satisfying $g(S,S)=-g(T,T)=\frac12$
such that $V=S+T$. Define $U=S-T$, which is a null vector such that $g(U,V)=g(S,S)-g(T,T)=1$, and consider a pseudo-orthonormal
basis $\{U,V,E_1,\dots,E_{n}\}$. For any vector field $Z$, from equations \eqref{curv-simplify} and \eqref{curv-general} we get
{
\begin{align}\label{eq:formula-curv}
R(Z,E_i,E_j,V) & = -\frac{1}{n+1}\rho(Z,V)\delta_{ij}+\frac{1}{n+1}\rho(E_i,V)g(Z,E_j)\\
    &  = \frac{\tau}{n(n+1)}g(Z,V)\delta_{ij}-\frac{\tau}{n(n+1)}g(E_i,V)g(Z,E_j)\nonumber\\
    & \quad -\frac{1}{n}\rho(Z,V)\delta_{ij}-\frac{1}{n}\rho(E_i,E_j)g(Z,V)\nonumber\\
    & \quad +\frac{1}{n}\rho(Z,E_j)g(E_i,V)+\frac{1}{n}\rho(E_i,V)g(Z,E_j).\nonumber
\end{align}
}

We use the fact that $\operatorname{Ric}(V)=\lambda V$ (Remark~\ref{remark:eigenvectors}) to see that
\[
\rho(V,V)=0,\,\rho(U,V)=\lambda,\, \rho(V,E_i)=0 \text{ for all } i=1,\dots,n.
\]
On the other hand compute $R(U,E_i,E_j,V)$ in expression \eqref{eq:formula-curv} to get that

\begin{align*}
R(U,E_i,E_j,V)  =& -\frac{1}{n+1}\lambda\delta_{ij}\\
     =& \frac{\tau}{n(n+1)}\delta_{ij}-\frac{1}{n}\lambda\delta_{ij}-\frac{1}{n}\rho(E_i,E_j).
\end{align*}
Hence $\rho(E_i,E_j)=0$ if $i\neq j$ and $\rho(E_i,E_i)=\frac{\tau-\lambda}{n+1}$ for all $i=1,\dots, n$. Now, compute
\[
\tau=2\rho(U,V)+n\rho(E_i,E_i)=\frac{(n+2)\lambda+n\tau}{n+1}
\]
to see that $\tau=(n+2)\lambda$. Hence the scalar curvature $\tau$ is constant and from \eqref{condgradsol1} we have $0=\nabla \tau=2\operatorname{Ric}(V)=2\lambda V$. Therefore we conclude that $\lambda=0=\tau$ and the only possibly nonzero Ricci component is $\rho(U,U)$, so $(M,g,f)$ is a steady gradient Ricci soliton with nilpotent Ricci operator.

Since the soliton is steady, from (\ref{gradsoliton}) we have
{ $\operatorname{hes}_f=-\operatorname{Ric}$.}
Now since $\operatorname{Ric}(V)=0$, it follows that $\nabla_{V}V=0$, which shows that $V$ is a geodesic vector field.

The gradient of the potential function is a recurrent vector field (i.e., the null line field $\mathcal{D}=\text{span}\{\nabla f\}$ is parallel) if and only if { $\nabla_X \nabla f=\operatorname{hes}_f(X)=\sigma(X)\nabla f$ for some $1$-form $\sigma$ and for all $X$}. Since $(M,g,f)$ is a steady gradient Ricci soliton, it follows from the expressions above for the Ricci operator that
\begin{align*}
&\text{hes}_f(U)=-\operatorname{Ric}(U)=-\rho(U,U) V,\\
&\text{hes}_f(V)=-\operatorname{Ric}(V)=0,\\
&\text{hes}_f(E_i)=-\operatorname{Ric}(E_i)=0,
\end{align*}
showing that $V$ is a recurrent vector field with $1$-form { $\sigma$ given by $\sigma(U)=-\rho(U,U)$, $\sigma(V)=0$ and $\sigma(E_i)=0$} for all $i=1,\dots,n$.

It follows now from (\ref{curv-general}), the expressions of the Ricci tensor above and the vanishing of the scalar curvature that
\[
R(\mathcal{D}^\perp,\mathcal{D}^\perp,\cdot,\cdot)=0.
\]
This shows that $(M,g)$ is a $pr$-wave. Moreover note that the Ricci tensor is isotropic and thus that $(M,g)$ is indeed a $pp$-wave
\cite{leistner}.
\end{proof}

\begin{remark}\rm
Note that although $(M,g)$ is a $pp$-wave, and hence it admits a null parallel vector field, $\nabla f$ is not in general parallel.
\end{remark}

\section{Ricci solitons on $pp$-waves}\label{se:5}

In this section we analyze the existence of gradient Ricci solitons on $pp$-waves. { Theorem \ref{mainth}-$(ii)$ will follow as a consequence of Lemma~\ref{lemma:isotropic} and the analysis performed here.} $pp$-waves are the underlying structure corresponding to many Lorentzian geometric properties without Riemannian analog, thus they are a natural family to look for new examples of complete gradient Ricci solitons.

{ Henceforth} we set $M=\mathbb{R}^{n+2}$ with coordinates $(u,v,x_1,\dots,x_{n})$, and metric { $g_{ppw}$} given by \eqref{pp} for some arbitrary function $H(u,x_1,\dots,x_{n})$. The possibly nonzero components of the Levi-Civita connection in the basis of coordinate vector fields $\{\partial_u=\frac{\partial}{\partial _{u}}, \partial_v=\frac{\partial}{\partial _{v}}, \partial_i=\frac{\partial}{\partial _{x_i}}\}$ are
\begin{equation}\label{nablaei}
\nabla _{\partial_u} \partial _u= \frac{1}{2}\partial_u H\, \partial _{v}-\frac{1}{2}\sum _{i=1}^n \partial_i H\, \partial _{i},  \qquad
\nabla _{\partial_u} \partial _i= \frac{1}{2}\partial_i H\, \partial _{v},  \quad i=1,\dots,n.
\end{equation}
This shows that the { null} vector field $\partial_v$ is parallel. Moreover, the possibly nonvanishing components of the curvature tensor are given (up to the usual symmetries) by
\begin{equation}\label{curvcomp}
R_{uiuj}=-\frac 12 \partial_{ij}^2 H,  \quad i,j=1,\dots,n.
\end{equation}
The scalar curvature $\tau$ is zero, since the Ricci tensor is determined by
 \begin{equation}\label{rhoei}
\rho _{uu} = -\frac{1}{2} \sum_{i=1}^n \partial_{ii}^2 H.
\end{equation}
Therefore, a $pp$-wave is Einstein (and hence Ricci flat) if and only if the space-Laplacian of the defining function $H$ vanishes identically.

\begin{theorem}
{ $(M,g_{ppw},f)$} is a nontrivial gradient Ricci soliton if and only if it is steady and the potential function $f$ satisfies  $f(u,x_1,\dots,x_n)=f_0(u)+\underset{i=1}{\overset{n}{\sum }} \kappa_i x_i$, where { $\kappa_i$ are arbitrary constants and}
\[
f_0''(u)=-\rho_{uu}-\frac 12 \sum_{i=1}^n \kappa_i \partial_i H(u,x_1,\dots,x_n).
\]
\end{theorem}

\proof
Let $f$ be a function on $\mathbb{R}^{n+2}$.
Then the gradient is given by
$\nabla f=(\partial_vf,\partial_u f-H \partial_vf,\partial_1f,\dots,\partial_nf)$
and thus (\ref{gradsoliton}) becomes
{ \begin{equation}\label{sysgradsol}
\left\{ \begin{array}{lc}
\frac 12\underset{i=1}{\overset{n}{\sum }} \partial_iH\, \partial_if+\partial^2_{uu}{f}-\frac 12\partial_{u}H\, \partial_v f+\rho_{uu}=\lambda H, &\quad  \\
\noalign{\bigskip}
\partial^2_{ui}f-\frac 12\partial_iH\, \partial_v f=0, &\quad 1\leq i\leq n,\\
\noalign{\bigskip}
\partial_{ii}^2f=\lambda, &\quad 1\leq i\leq n,\\
\noalign{\bigskip}
\partial_{uv}^2f=\lambda, &\quad \\
\noalign{\bigskip}
\partial^2_{ij}f=\partial^2_{vi}f=\partial^2_{vv}f=0, &\quad 1\leq i\neq j\leq n.
\end{array}\right.
\end{equation}}
Integrating {  equations} $\partial^2_{vi}f=\partial^2_{vv}f=0$ in \eqref{sysgradsol} we obtain that the potential function splits as $f(u,v,x_1,\dots,x_n)=f_0(u,x_1,\dots,x_n)+v f_1(u)$ for some functions $f_0$, $f_1$.
Moreover equations $\partial^2_{uv}f=\lambda$ and $\partial^2_{ij}f=0$ now show that
$f(u,v,x_1,\dots,x_n)=\underset{i=1}{\overset{n}{\sum }} f_i(u,x_i) +v (\lambda u+\kappa)$ for some constant $\kappa$
and functions $f_i$, $i=1,\dots,n$. Hence \eqref{sysgradsol} reduces to
\begin{equation}\label{sysgradsol2}
\left\{
\begin{array}{lc}
\frac 12\underset{i=1}{\overset{n}{\sum }}\partial_iH\,\partial_if_i+\underset{i=1}{\overset{n}{\sum }}\partial^2_{uu}{f_i}-\frac 12(\lambda u+\kappa)\partial_{u}H+\rho_{uu}=\lambda H, &\\
\noalign{\bigskip}
\partial^2_{ui}f_i-\frac 12(\lambda u+\kappa)\partial_iH=0, & \quad 1\leq i\leq n, \\
\noalign{\bigskip}
\partial^2_{ii}f_i=\lambda, & \quad 1\leq i\leq n.
\end{array}
\right.
\end{equation}
Integrating the last equations in \eqref{sysgradsol2} { we have}
\[
f_i(u,x_i)=f_{0,i}(u)+x_i \kappa_i(u)+\frac{\lambda}{2} x_i^2,
\]
for some functions $f_{0,i}(u)$ and $\kappa_i(u)$.
Substituting the above into \eqref{sysgradsol2} and differentiating the second set of equations one gets
\[
0=\partial^3_{uii}f_i=(\lambda u+\kappa)\partial_{ii}^2 H,
\]
which shows that either,  $\partial_{ii}^2 H=0$ for all $i$ (and hence the $pp$-wave is Ricci flat) or otherwise $\lambda=\kappa=0$.

The first case, when $(M,g_{ppw})$ is Ricci flat, was already analyzed in Theorem~\ref{th:Einstein}.
The second case, $\lambda=\kappa=0$ shows that non Einstein gradient Ricci solitons are steady and $f$ becomes $f(u,v,x_1,\dots,x_n)=f_0(u)+\underset{i=1}{\overset{n}{\sum }} \kappa_i(u) x_i$. Now the second equation in \eqref{sysgradsol2} reduces to
$\kappa_i^\prime(u) = 0$
and hence $\kappa_i(u)=\kappa_i$ for  real constants $\kappa_i$, which gives
\begin{equation}\label{eq:f}
f(u,v,x_1,\dots,x_n)=f_0(u)+\underset{i=1}{\overset{n}{\sum }} \kappa_i x_i.
\end{equation}

Finally, it follows from the first equations in \eqref{sysgradsol2} that the function $f_0(u)$ is given by the differential equation
\begin{equation}\label{condult}
f_0^{\prime \prime}=\frac 12 \left(\sum_i \partial^2_{ii}H\right)-\frac 12\sum_i \kappa_i \partial_iH
=-\rho_{uu}-\frac 12 \sum_i \kappa_i \partial_iH,
\end{equation}
which completes the proof.
\qed

\begin{remark}\rm
In general, equation (\ref{condult}) does not have a solution, since the derivatives $\partial_i H(u,x_1,\dots,x_n)$ and $\partial_{ii}H$ may be functions of the $x_i$'s. Further note that $\nabla f$ is not isotropic in general since
$\|\nabla f\|=\underset{i=1}{\overset{n}{\sum }} \kappa_i^2$, although it is a geodesic vector field
since $\nabla_{\nabla f}\nabla f=-\operatorname{Ric}(\nabla f)=0$.
\end{remark}

\begin{remark}\rm
A special class of $pp$-waves are \emph{plane waves}, which are defined by setting
\begin{equation}\label{plane-pp}
H(u,x_1,..,x_{n})=\sum_{i,j} a_{ij}(u)x_ix_j.
\end{equation}
Note that any plane wave is a isotropic steady gradient Ricci soliton for a potential function $f$  given by \eqref{eq:f} for constants $\kappa_i=0$ for all $i=1,\dots, n$, since \eqref{condult} becomes
\[
f_0^{\prime \prime}(u)=\frac 12 \left(\sum_i a_{ii}(u)\right).
\]
Moreover, it is shown in \cite{CFS} that plane waves are geodesically complete and therefore since $\nabla f$ is a geodesic vector field, it follows that $\nabla f$ is complete.
\end{remark}

\subsection{Locally conformally flat $pp$-waves}

It follows from the expressions { \eqref{curvcomp} and \eqref{rhoei}}, that a $pp$-wave is locally conformally flat if and only if the defining function $H$ takes the form
\begin{equation}\label{cflat}
H(u,x_1,\dots,x_n)= a(u) \sum _{i=1}^n x_i ^2 + \sum _{i=1}^n b_i (u) x_i +c(u),
\end{equation}
where $a,b_1,\dots,b_n,c$ are smooth functions of one variable.

In this case condition (\ref{condult}) reduces to
\begin{eqnarray}\label{condult-1}
f_0^{\prime \prime}=-\rho_{uu}-\frac 12 \sum_{i=1}^n \kappa_i b_i(u)- a(u)\sum_{i=1}^n \kappa_i x_i,
\end{eqnarray}
where $\rho_{uu}=-n a(u)$. So, if we differentiate \eqref{condult-1} with respect to $x_i$ we get that $a(u) \kappa_i=0$ for all $i=1,\dots,n$. Then, unless the manifold is flat, it follows that  necessarily $\kappa_i=0$ for all $i$, and the potential function is given by
\[
f(u,v,x_1,\dots,x_n)=f_0(u), \quad \mbox{where}\quad  f_0''(u)=-\rho_{uu}=n a(u).
\]
{ This completes the proof of Theorem~\ref{mainth}.}

\begin{remark}\rm
Gradient Ricci solitons are a particular family of Ricci solitons, which are triples $(M,g,X)$ where $X$ is a vector field satisfying:
\begin{equation}\label{soliton}
\frac 12{\mathcal L}_{X}g+\rho=\lambda g.
\end{equation}
Here ${\mathcal L}_{X}$ denotes the Lie derivative in the direction of $X$. Note that when $X=\nabla f$ equation \eqref{soliton} becomes equation \eqref{gradsoliton}.

{
Let $\displaystyle X=X_u \partial_u+X_v \partial_v+\sum_i X_i \partial_{i}$ be an arbitrary vector field on $(\mathbb{R}^{n+2},g_{ppw})$. Then
\eqref{soliton} becomes
\begin{equation}\label{syssolconfflat}
\left\{\begin{array}{lc}
\frac 12\underset{i=1}{\overset{n}{\sum }} \partial_iH X_i+\frac 12\partial_{u}H X_u+H \partial_u X_u+\partial_u X_v+\rho_{uu}=\lambda H,\!&\\
\noalign{\bigskip}
H \partial_v X_u+\partial_v X_v+\partial_u X_u=2\lambda, \!& \\
\noalign{\bigskip}
H \partial_i X_u+\partial_i X_v+\partial_u X_i=0, \!&\! 1\leq i\leq n,\\
\noalign{\bigskip}
\partial_i X_j+\partial_j X_i=0, \!&\! 1\leq i\neq j \leq n, \\
\noalign{\bigskip}
\partial_i X_u+\partial_v X_i=0, \!&\! 1\leq i \leq n, \\
\noalign{\bigskip}
\partial_v X_u=0; \quad \partial_i X_i=\lambda, \!&\! 1\leq i\leq n.
\end{array}\right.
\end{equation}}

Consider the vector field
\begin{equation}\label{eq:part-solution-confflat}
X=\left(p(u)-\underset{i=1}{\overset{n}{\sum }} q_i'(u) x_i+2\lambda v\right) \partial_{v}
+
\underset{i=1}{\overset{n}{\sum }} \left(q_i(u)+\lambda x_{i}\right) \partial_{i},
\end{equation}
where functions $p$ and $q_i$ satisfy the following conditions
\begin{equation}\label{functcond}
\left\{\begin{array}{lc}
a(u) q_i(u)-q_i''(u)=\frac{\lambda}{2}b_i(u), &\quad 1\leq i \leq n,\\
\noalign{\bigskip}
\frac 12 \underset{i=1}{\overset{n}{\sum }}b_i(u) q_i(u)+\rho_{uu}+p'(u)=\lambda c(u). &
\end{array}\right.
\end{equation}
{ Note that one can always find $p$, $q_i$ being solutions of \eqref{functcond}. A straightforward calculation from \eqref{syssolconfflat} shows that $(M,g,X)$ is a Ricci soliton. Also observe} that $\lambda$ is the constant of equation \eqref{soliton} and can be chosen with absolute freedom.
In contrast to the gradient case, we obtain that
\emph{any locally conformally flat $pp$-wave $(M,g)$ admits appropriate vector fields resulting in expanding, steady and shrinking Ricci solitons}.
\end{remark}

\subsection{Lorentzian manifolds with recurrent curvature}

A pseudo-Riemannian manifold $(M,g)$ is said to be \emph{recurrent} (or with \emph{recurrent curvature}) if $\nabla R=\sigma\otimes R$ for some $1$-form $\sigma$.

Observe here that \emph{locally conformally flat $pp$-waves are recurrent}. Indeed it suffices to show that the Ricci tensor is recurrent due to local conformal flatness. For any locally conformally flat $pp$-wave $(M,g_{ppw})$ as in \eqref{cflat} the only nonzero component of the Ricci tensor is $\rho_{uu}=-na(u)$, and hence from \eqref{nablaei} it follows that the only nonzero component of $\nabla\rho$ becomes
$\nabla_{\partial_u}\rho_{uu}=-na'(u)$. Hence $\nabla\rho=\sigma\otimes\rho$, just considering the
$1$-form $\sigma=(\ln{a(u)})'\, du$.

Recurrent Lorentz manifolds have been classified by Walker \cite{walker-recurrent} (see also Galaev  \cite{galaev-recurrent}). Non-symmetric Lorentzian recurrent manifolds are $pp$-waves which correspond to one of the following two families
\begin{description}
\item [Type I] The defining function satisfies $H(u,x_1,\dots,x_n)=H(u,x_1)$ where $\partial_{11}^2 H$ is not constant.
\item [Type II] The defining function is given by
$H(u,x_1,\dots,x_n)=a(u)\left(\sum_{i=1}^n b_i x_i^2\right)$
for constants $b_1,\dots,b_n$ with $|b_1|\geq\dots\geq |b_n|,\, b_2\neq 0,$ and a function $h$ such that
       $a'(u)\neq 0$.
\end{description}

For a recurrent manifold of Type I, condition (\ref{condult}) reduces to
{ \begin{eqnarray}\label{condult-2}
f_0^{\prime \prime}(u)=-\rho_{uu}-\frac 12 \kappa_1 \partial_1 H(u,x_1),
\end{eqnarray}}
where $\rho_{uu}=-\frac 12 \partial_{11}^2 H(u,x_1)$.

{ Differentiating} in \eqref{condult-2} with respect to $x_1$ we get  $\frac{\kappa_1}{2}\partial_{11}^2 H(u,x_1)$ $-$ $\frac 12 \partial_{111}^3 H(u,x_1)=0$, and hence
the defining function $H(u,x_1)$ becomes
\[
H(u,x_1)=\frac{1}{\kappa_1^2}e^{\kappa_1 x_1} h_0(u)+h_1(u)+x_1 h_2(u).
\]
 Then the soliton is given by (see \eqref{eq:f}-\eqref{condult}) $f(u,v,x_1,\dots,x_n)=f_0(u)+\underset{i=1}{\overset{n}{\sum }} \kappa_i x_i,$ where
\[
f_0''(u)=-\frac{\kappa_1}{2} h_2(u).
\]

Note that in this case { $\nabla f$} is always spacelike (since $\kappa_1\neq 0$) and the underlying
manifold is not locally conformally flat (unless it is flat which occurs if $h_0(u)= 0$).

\smallskip

For a recurrent manifold of Type II condition (\ref{condult}) reduces to
\begin{eqnarray}\label{condult-3}
f_0^{\prime \prime}(u)=-\rho_{uu}-a(u)\sum_{i=1}^n \kappa_i b_i x_i
\end{eqnarray}
where $\rho_{uu}=-a(u)\sum_i b_i$.
Taking derivatives with respect to $x_i$ in \eqref{condult-3} we get that $\kappa_i b_i a(u)=0$ for all $i$ and therefore, as $b_1\neq 0 \neq b_2$ there are two different possibilities:
\begin{itemize}
\item[(a)] If { $b_i\neq 0$  for all $i$ then,} unless the manifold is flat, it follows that $\kappa_i=0$ for all $i$ and the potential function is given by $f(u,v,x_1,\dots,x_n)=f_0(u)$ where $f_0''(u)=-\rho_{uu}=a(u)\sum_{i=1}^nb_i$. In this case the Ricci soliton is isotropic ($\|\nabla f\|=0$).

\item[(b)] { If $b_j=0$} for some $j\in \{3,\dots,n\}$,
then $\kappa_i=0$ for $i< j$ and the potential function is given by
\[
f(u,v,x_1,\dots,x_n)=f_0(u)+\underset{i=j}{\overset{n}{\sum }} \kappa_i x_i\,,\mbox{ where }
f_0''(u)=-\rho_{uu}=a(u)\sum_{i=1}^jb_i.
\]
Further observe that in this case $\nabla f$ is spacelike.
\end{itemize}

Summarizing the above, we have that
\begin{description}
\item[Type I]
A recurrent Lorentzian manifold of Type I admits a function $f$ resulting in a gradient Ricci soliton if and only if
the defining function $H(u,x_1)$ satisfies
$H(u,x_1)=\frac{1}{\kappa_1^2}e^{\kappa_1 x_1} h_0(u)+h_1(u)+x_1 h_2(u)$. Moreover, in this case
{ $\nabla f$} is a spacelike vector field.

\item[Type II]
A recurrent Lorentzian manifold of Type II always admits  a function $f$ resulting in a gradient Ricci soliton. The causal character of $\nabla f$ may be null or spacelike.
\end{description}

%

\subsection{Two-symmetric Lorentzian manifols}

As a generalization of locally symmetric spaces, Lorentzian manifolds whose higher order derivatives of the curvature tensor vanish, have been investigated. A Lorentzian manifold is said to be  \emph{two-symmetric} if $\nabla^2 R=0$ but $\nabla R\neq 0$. It was shown by Senovilla \cite{senovilla} that two-symmetric Lorentzian manifolds admit a parallel null vector field and the local structure of such manifolds was given recently in \cite{BSS}, \cite{alekseevsky-galaev} showing that they are a special family of $pp$-waves.

An $(n+2)$-dimensional Lorentzian manifold is two-symmetric if and only if it is a $pp$-wave as in \eqref{pp} with
\begin{equation}\label{g2sym}
        H(u,x_1,\dots,x_n)=\sum_{i,j=1}^n (a_{ij}u+b_{ij})x_ix_j\,,
\end{equation}
where $(a_{ij})$ is a diagonal matrix with the diagonal elements $a_{11}\leq\dots \leq a_{nn}$ non-zero real numbers and $(b_{ij})$ an arbitrary symmetric matrix of real numbers.

Now, an immediate application of \eqref{eq:f}-\eqref{condult} shows that two-symmetric Lorentzian manifolds are gradient Ricci solitons whose potential function is given by
$f=f_0(u)$ where
\[
 f_0''(u)=-\rho_{uu}=\underset{i=1}{\overset{n}{\sum }} (b_{ii}+u a_{ii}).
\]
Finally observe that $\nabla f$ is a geodesic vector field and $(M,g)$ is geodesically complete. Moreover $\nabla f$ is isotropic.

\subsection{Conformally symmetric Lorentzian manifolds}

A Lorentzian manifold is said to be \emph{conformally symmetric} if the covariant derivative of the Weyl tensor vanishes identically ($\nabla W=0$). Clearly locally symmetric and locally conformally flat manifolds are conformally symmetric and the converse is true in the Riemannian setting. In Lorentzian signature there exist, however, conformally symmetric manifolds which are neither locally conformally flat nor locally symmetric. These manifolds have recurrent Ricci tensor and have been described locally by Derdzinski and Roter \cite{DR1}. It turns out that all of them are $pp$-waves given by
\[
H(u,x_1,\dots,x_n)=a(u)\sum_i x_i^2+\sum_{i,j}b_{ij} x_ix_j
\]
where $(b_{ij})$ is a nonzero symmetric matrix with { $\sum_{i=1}^nb_{ii}=0$}.

Now equations \eqref{eq:f} and \eqref{condult} show that any conformally symmetric Lorentzian manifold of this family admits a function $f$ resulting in a steady gradient Ricci soliton with isotropic $\nabla f$.


\end{document}